\documentclass[3p,times]{elsarticle}
\usepackage{dsfont}
\usepackage{amsmath, mathtools}
\usepackage{amssymb}
\usepackage{graphicx}
\usepackage{enumitem} 
\usepackage{hyperref}
\usepackage{cleveref}
\usepackage{amsthm}
\usepackage{enumerate}
\usepackage{multirow}
\usepackage{algorithm}
\usepackage{algpseudocode}
\usepackage{xcolor}
\theoremstyle{definition}
\newtheorem{Theorem}{Theorem}[section]
\newtheorem{Proposition}{Proposition}[section]%[theorem]

\newtheorem{Example}{Example}

\usepackage[all]{xy}
\newtheorem{Definition}{Definition}[section]
\newcommand{\norm}[1]{\left\lVert#1\right\rVert}
\journal{}

\newenvironment{Proof}{\paragraph{Proof:}}{\hfill$\blacksquare$}
\newenvironment{prof}{\paragraph{Sketch of proof for pointwise convergence:}}{\hfill$\blacksquare$}

\newcommand\proba{{\operatorname{\mathbb{P}}}}
\newcommand\E{{\operatorname{\mathbb{E}}}}
\newcommand\var{{\operatorname{\mathrm{Var}}}}
\newcommand{\ZCRR}{Z_{n}^{\mathrm{CRR}}}

\newcommand{\Csubbin}{C_{{S}}^\textrm{CRR}}
\newcommand{\Csubbs}{C_{{S}}}
\newcommand{\floor}[1]{\lfloor #1 \rfloor}

\usepackage{setspace}
\begin{document}

\begin{frontmatter}

\title{About subordinated generalizations of 3 classical models of option pricing}
\author[label1]{Micha\l{} Balcerek}
\ead{michal.balcerek@pwr.edu.pl}

\author[label1]{Grzegorz Krzy\.zanowski}
\ead{grzegorz.krzyzanowski@pwr.edu.pl}

\author[label1]{Marcin Magdziarz}
\ead{marcin.magdziarz@pwr.wroc.pl}

\address[label1]{Hugo Steinhaus Center,
Faculty of Pure and Applied Mathematics, Wroclaw University of Science and Technology
50-370 Wroclaw, Poland}

\begin{abstract}
    In this paper, we investigate the relation between  Bachelier and Black-Scholes models driven by the infinitely divisible inverse subordinators. Such models, in contrast to their classical equivalents, can be used in markets where periods of stagnation are observed. We introduce the subordinated Cox-Ross-Rubinstein model and prove that the price of the underlying in that model converges in distribution and in Skorokhod space to the price of underlying in the subordinated Black-Scholes model defined in \cite{gajda}. Motivated by this fact we price the selected option contracts using the binomial trees. The results are compared to other numerical methods.
\end{abstract}

\begin{keyword}
%% keywords here, in the form: keyword \sep keyword
subdiffusion, subordinator, option contracts, Monte Carlo, Bachelier model, Black-Scholes model, Cox-Ross-Rubinstein model 

%% MSC codes here, in the form: \MSC code \sep code
%% or \MSC[2008] code \sep code (2000 is the default)

\end{keyword}

\end{frontmatter}

% Use "Eq" instead of "Equation" for equation citations.
\section*{Introduction}
\label{intro}
The theories of finance and economy began nearly a century ago with the works of Walras and Pareto and their mathematical description of equilibrium theory (see \cite{focardi} and references therein). Soon after L. Bachelier in his thesis \cite{Luis} began the theory of option pricing. 
In the same work, he initiated the research of diffusion processes five years before recognized as the pioneering works of A. Einstein \cite{Einstein}, M. Smoluchowski \cite{Smoluchowski} and decades before famous works of P. L\'{e}vy \cite{Levy}, K. It\^{o} \cite{Ito} and N. Wiener \cite{Wiener}.
Bachelier was the first who found the practical application of random walk in finance.
Despite the groundbreaking character, the theory of Bachelier was forgotten for over fifty years. In $1965$ Samuelson \cite{samuelson2015rational} used the Geometric Brownian Motion (GBM) to model prices of underlying instrument.
The culminating event in the developmenet of the theory of option pricing was $1973$ when Black, Scholes and Merton found consistent formulas for the fair prices of European options \cite{black, Merton}.
The discovery was of such great importance that Merton and Scholes were awarded the Nobel Prize in Economics in $1997$.

The Bachelier model is intuitive and elegant, but it allows the underlying financial process to drop into the negative regime. The underlying asset in B-S model instead of the Arithmetic Brownian Motion (ABM) follows GBM. Hence, the obvious property of non-negativity 
of the underlying process is conserved. %Moreover, the classical Bachelier model assumes that interest rate $r=0$ \cite{musiela}.
The classical Black-Scholes (B-S) model, despite its popularity, can not be used in many different cases. Note that one of the greatest weaknesses of B-S model is the presence of volatility smile - in other words the assumption of constant volatility results that B-S price often does not fit the market price \cite{hull2003options}. Therefore, the B-S model was generalized to weaken its strict assumptions, allowing such features as
the possibility of jumps \cite{chevallier2017estimation, zhen2020dissecting}, stochastic volatility \cite{alghalith2020pricing,  zhen2020dissecting}, market regulations \cite{aliahmadi2020option},  transactions costs \cite{liu2013closed, tian2020european}, regime switching \cite{lin2020regime} and stagnation of underlying asset \cite{ja3, ja2, MM} to name only few.

Nowadays, we know that financial data are often far away from classical models. Neither the B-S nor the Bachelier model can capture the characteristic motionless periods observed among many markets (see, e.g., \cite{AW11} and the references therein). Such dynamics can be observed in every nonliquid market, e.g., in emerging markets where the number of transactions and participants is low. Similar behaviour can be observed in experimental physics where the motion of small particles is interrupted by trapping events. During that random time, the particle is immobilized and stays motionless. For the discussed systems, the subdiffusive model was proposed and nowadays it is well-established approach with many practical applications \cite{AW17,AW16, AW15, AW14, AW18}.

In order to properly price the option contracts in markets where the subdiffusive character of the underlying asset was observed, the generalization of Bachelier \cite{orzel} and B-S model \cite{MM} were proposed. 
For the subdiffusive B-S model, effective numerical methods for valuing vanilla and barrier options were proposed \cite{ja3, ja2}.
The subdiffusive models of Bachelier and B-S are generalizations of the standard models to the case where the underlying assets display characteristic periods in which they stay motionless. The classical models do not take these phenomena into account. As a consequence of an option valuation for the market displaying the stagnation property, the fair price provided by the Bachelier/B-S model is misestimated.
To describe these dynamics properly, the subordinated Bachelier/B-S model assumes that the underlying instrument is driven by an inverse subordinator.

Let us focus on the inverse tempered stable subordinator. In this case, time periods in which their value are observed to be constant depend on the subdiffusion parameter $\alpha\in(0,1)$ and tempering parameter $\lambda\geq0$. If $\alpha=1$, the subordinated Bachelier/B-S reduce to the classical models \cite{MM, orzel}.

\begin{figure}[ht]
	%\centering
	\raggedleft
	\includegraphics[scale=0.38]{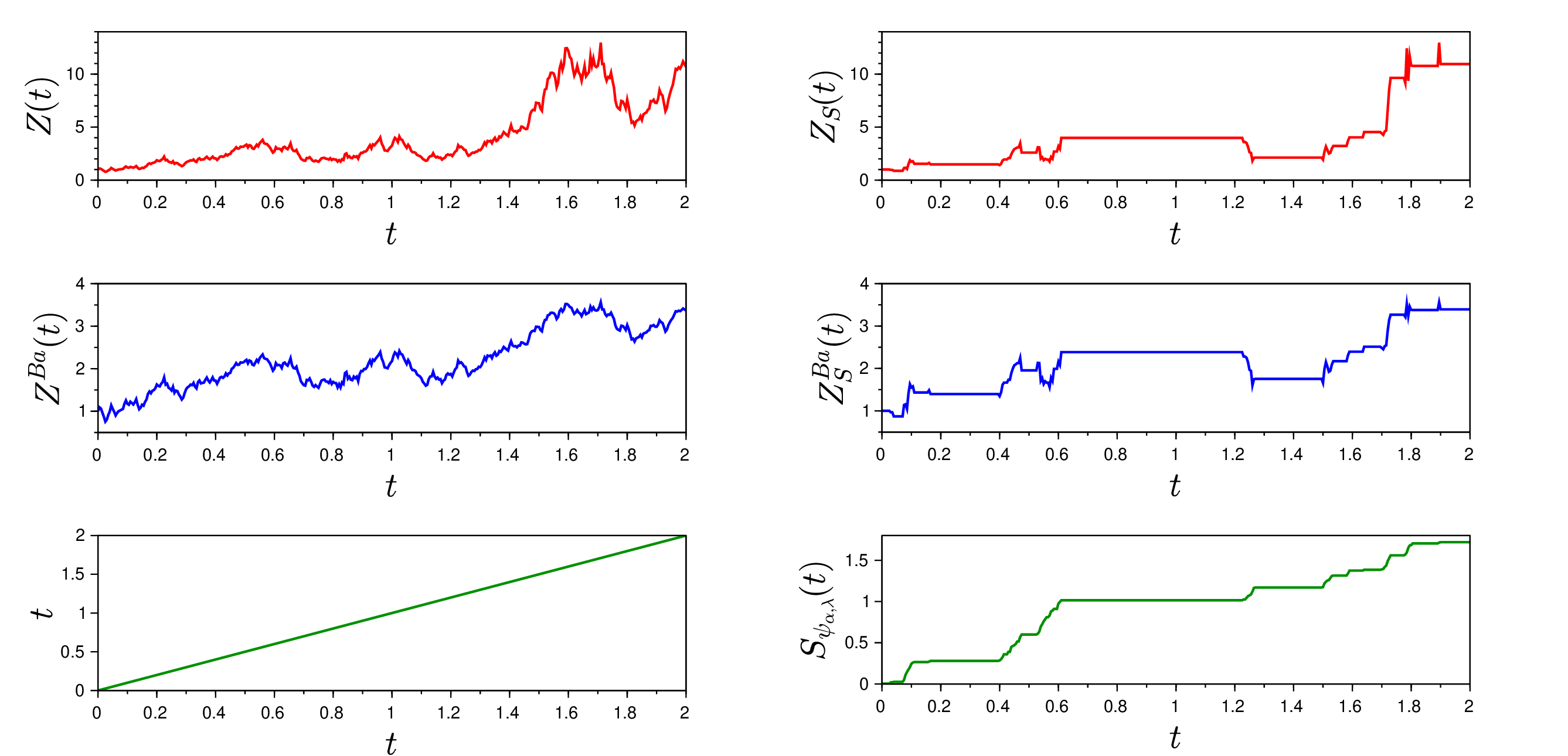}
	\caption{{\bf The sample path of GBM (top) and ABM (middle) for the same realization of $B(t)$ (left)/$B(t)$ and $S_{\psi_{\alpha, \lambda}}(t)$ (right).} The classical models (left), in contrast to their subordinated equivalents (right), can not predict even short stagnation period. 
		Each flat period of inverse subordinator (right bottom) induces flat period in GBM/ABM. The parameters are $Z_{0}=\mu=\sigma=\lambda=1$, $\alpha=0.7$.}  
	\label{fig:paths}
\end{figure}

In Figure \ref{fig:paths} we compare sample trajectories of the underlying assets in the classical Bachelier and B-S models compared to their equivalents characterized by $\lambda$-tempered $\alpha$-stable inverse subordinator. Even a small stagnation of an underlying can not be simulated by a classical model.
As a generalization of the standard Bachelier and B-S models, their subordinated equivalents can be applied in a wide range of markets - including all cases where Bachelier and B-S can be used. 

The Cox-Ross-Rubinstein (CRR) model is a discrete version of B-S model and provides a general and easy way to price a wide range of option contracts \cite{cox1979option}. 

In this paper, we find the relation between the European options prices in the subordinated Bachelier and B-S model.
By this way, we show that the earliest and the most popular models of option pricing are in some sense close even if they are run by inverse subordinators.
We also define the subordinated CRR (called also subordinated binomial) model and prove that the price of the underlying/European option converges to the price of the underlying/European option in the subordinated B-S model. We use this fact to price the selected types of options. Finally, we compare the results with other numerical methods.

\section{Inverse subordinators}
We begin with defining the Laplace exponent $\displaystyle\psi(t)$
\begin{equation}
\displaystyle\psi(t)=\int_{0}^{\infty}(1-e^{-ux})\nu(dx),\label{tfutfu}
\end{equation} 
where $\nu$ is the L\'{e}vy measure such that 
$\displaystyle\int_{0}^{\infty}(\min(1,x))\nu(dx)~<~\infty$. We assume that $\displaystyle\nu(0,\infty)=\infty$.
Then the so-called inverse subordinator $S_{\psi}(t)$ can be defined as  \cite{magdziarz2009langevin, gajda, orzel}
\begin{equation}
	\displaystyle S_{\psi}(t)=\inf\{\tau>0: U_{\psi}(\tau)>t\},\label{Def_subordynatora}
\end{equation}
where $U_{\psi}(t)$ is the strictly increasing L\'{e}vy process (subordinator) \cite{ken1999levy} with the Laplace transform 
\begin{equation*}
	\displaystyle\E \exp(-u U_{\psi}(t))=\exp(-t\psi(u)),
\end{equation*}
and $\psi(u)$ is given by (\ref{tfutfu}). Since $\displaystyle U_{\psi}(t)$ is a pure jump process with c\'adl\'ag trajectories and not compound Poisson, $S_{\psi}(t)$ has continuous and singular sample paths with respect to the Lebesgue measure. Moreover, every jump of $U_{\psi}(t)$ correspond to the flat period of inverse. Note that inverse subordinators except of compound Poisson (which does not satisfy the assumption $\displaystyle\nu(0,\infty)=\infty$) have density \cite{meerschaert2008triangular}.
In the next sections we use the notation $S(t)=S_{\psi}(t)$. In the whole paper, we call the model ``subdiffusive'' if it is driven by the $\alpha$-stable inverse subordinator, i.e. if $\psi(u)=u^{\alpha}$, $\alpha\in(0,1)$. We call the model ``subordinated'' if it is driven by any inverse subordinator $S(t)$. Note that in the literature subordinated models are the models run by the subordinators or, like in this paper, by the inverse subordinators.

\section{Subordinated Bachelier model}
Bachelier's model is based on the observation that the stock price dynamics is analogous to the motion of small particles suspended in liquids \cite{Luis}.
Therefore, in modern terms, the stock price follows ABM 
\begin{equation*}
	Z^{Ba}(t)=Z_{0}+\mu t+\sigma B(t),
\end{equation*}
where $Z^{Ba}(t)$ is the price of the underlying instrument, $t\geq0$, $Z_{0}=Z^{Ba}(0)$, $\mu\in(-\infty,\infty)$, $\sigma\in(0,\infty)$ and $B(t)$ is a standard Brownian motion on the probability space
$\left(\Omega, \mathcal{F} , \mathbb{P}\right)$.
Here, $\Omega$ is the sample space,  $\mathcal{F}$ is a $\sigma$-algebra interpreted as the information about the history of underlying price and $\mathbb{P}$ is the ``objective'' probability measure.

The price of an underlying asset in the subordinated Bachelier model is defined by the subordinated ABM: 
\begin{equation}
	\begin{cases}
		Z^{Ba}_{S}(t)=Z^{Ba}(S(t))=Z^{Ba}(0)+\mu S(t)+\sigma B(S(t)),\\
		Z^{Ba}\left(0\right)=Z_{0}, 
	\end{cases}\label{defsubABM1}
\end{equation} where $S(t)$ is defined in (\ref{Def_subordynatora}). We assume that $S(t)$ is independent  of $B(t)$ for each $t\in[0,T]$. Let us introduce the probability measure \begin{align}
	\mathbb{Q}^{Ba}\left(A\right)=\int_{A}\exp\left(-\gamma B\left(S\left(T\right)\right)-\frac{\gamma^{2}}{2}S\left(T\right)\right)d\mathbb{P},
	\label{1}
\end{align}
where $\gamma=\mu/\sigma$, $A\in\mathcal{F}$. As it is shown in \cite{orzel}, $Z^{Ba}_{S}(t)$ is a $\mathbb{Q}^{Ba}$-martingale.
The subordinated Bachelier model is arbitrage-free and incomplete \cite{orzel}. Let us recall that the lack of arbitrage is one of the most expected 
properties of the market and claims that there is no possibility to gain money without taking the risk. 
The market model is complete
if the set of possible gambles on the future states can
be constructed with existing assets. More formally, the market model is complete if every
$\mathcal{F}_{t\in[0,T]}$~-~measurable random variable $X$ admits a replicating self-financing strategy $\phi$ \cite{4}.
The Second Fundamental theorem of asset pricing \cite{4} states that a
market model described by $(\Omega, \mathcal{F}, \mathbb{P})$ and
underlying instrument $Z^{Ba}_{S}(t)$ with filtration $\mathcal{F}_{t\in[0,T]}$ is complete if and only if there is a unique
martingale measure equivalent to $\mathbb{P}$. 
Despite that $\mathbb{Q}^{Ba}$ defined in (\ref{1}) is not unique, it is the ``best'' martingale
measure in the sense of the criterion of minimal relative entropy. The relative entropy (we use the Kullback-Leibler approach) is the distance between measures \cite{kullback1951information}. The criterion of its minimal value in this case means, that the measure $\mathbb{Q}^{Ba}$ minimizes the distance to
the ``real-life'' measure $\mathbb{P}$ \cite{orzel}.

\section{Subordinated B-S model}\label{sectionbs}
The assumptions here are the same as in  the classical B-S model \cite{ja} with the exception that we do not assume the market liquidity and that 
the underlying instrument, instead of GBM, follows subordinated GBM \cite{gajda}:
\begin{equation}
	\left\{ \begin{array}{ll}
		Z_{S}\left(t\right)=Z\left(S(t)\right)=Z\left(0\right)\exp \left(\mu S(t)+\sigma B(S(t))\right),\\
		Z\left(0\right)=Z_{0}, 
	\end{array}\right.\label{defSubGBM}
\end{equation}
where $Z_{S}\left(t\right)$ - the price of the underlying instrument. The constants
$\mu$, $\sigma$ and the stochastic processes $S(t)$, $B(t)$ are defined in the same way as in the previous section. Note that $Z_{S}\left(t\right)$ equivalently can be defined by the subordinated ABM, i.e.,

$$Z_{S}\left(t\right)=Z_{0}\exp\left(Z^{Ba}_{S}(t)-Z_{0}\right),$$
where $Z^{Ba}_{S}(t)$ is defined in (\ref{defsubABM1}). Let us introduce the probability measure \begin{equation}
	\mathbb{Q}\left(A\right)=\int_{A}\exp\left(-\gamma B\left(S\left(T\right)\right)-\frac{\gamma^{2}}{2}S\left(T\right)\right)d\mathbb{P},\label{2}
\end{equation}
where $\gamma=(\mu-r+\frac{\sigma^2}{2})/\sigma$, $A\in\mathcal{F}$.
As shown in \cite{gajda}, $Z_{S}(t)$ is a $ \mathbb{Q}$-martingale.
The subordinated B-S model, similar to the subordinated Bachelier model, is arbitrage-free and incomplete \cite{gajda}. Despite that $ \mathbb{Q}$ defined in (\ref{2}) is not unique, it is the ``best'' martingale
measure in the sense of the criterion of minimal relative entropy \cite{gajda}. %It means that the measure $\mathbb{Q}$ minimizes the distance to
%the ``real life'' measure $\mathbb{P}$ \cite{gajda}. 
Note that for the case of stable and tempered stable inverse subordinators (where $S(t)=S_{{\psi}_{\alpha}}(t)$, $S(t)=S_{{\psi}_{\alpha,\lambda}}(t)$ respectively, for $\psi_{\alpha}(u)=u^{\alpha}$, $\psi_{\alpha,\lambda}(u)=(u+\lambda)^{\alpha}-\lambda^{\alpha}$, $\alpha\in(0,1)$, $\lambda\geq0$), for $\alpha=1$ $\mathbb{Q}^{Ba}$ and $\mathbb{Q}$ reduce to the measures of the classical Bachelier and B-S models which are arbitrage-free and complete \cite{MM, orzel}.
%\newpage

At the end of this section, we will find the relation between the European call and put options in the subordinated B-S model.
\begin{Proposition}
	Between the prices of the European call and put options in the subordinated B-S model, the following relation holds:
	\begin{equation}
		P_{S}(T)=C_{S}(T)+ K \E  e^{-rS(T)} - Z_{0},\label{putcallpar1}
	\end{equation}
	where $P_{S}(T)$ and $C_{S}(T)$ denote the fair prices of the European put and call options in the subordinated B-S model, respectively.
\end{Proposition}
\begin{Proof}
	By \cite{gajda}, we have
	\begin{equation}
		P_{S}(T)=\E P(S(T))\label{gajda11}
	\end{equation}
	
	and
	
	\begin{equation}
		C_{S}(T)=\E C(S(T)),\label{gajda22}
	\end{equation}
	where $C(T)$ and $P(T)$ denote the fair prices of the call and put options in the classical B-S model.
	By the put-call parity in the standard B-S model \cite{musiela}, for the given expiration time $\tau$ we have   \begin{equation}
		P(\tau)=C(\tau)+ K e^{-r\tau} - Z_{0}.\label{putcallpar2}
	\end{equation}
	After conditioning (\ref{putcallpar2}) on $\tau=S_{\psi}(T)$ and taking the expected value, by (\ref{gajda11}) and (\ref{gajda22}) the proof is completed. 
\end{Proof}

\section{Relation between the subordinated model of Bachelier and B-S}\label{section4}
Although the Bachelier and B-S models have a lot in common, they follow different concepts: an equilibrium argument in the Bachelier case, as opposed to the no-arbitrage argument in the B-S approach \cite{Sch}. The next theorem finds the direct relation between the European option prices in both models.\newpage

\begin{Theorem}
	Let us denote $T$ - expiration time, $\sigma$ - volatility of the B-S model, $C^{Ba}_{S}(T)$ and $C_{S}(T)$ - fair prices of a call option in the subordinated Bachelier and the subordinated B-S model with respect to the measures $\mathbb{Q}^{Ba}$ and $\mathbb{Q}$, respectively. We assume that $Z_{0}=K$ (at the money), $r=0$ and that $\sigma^{Ba}=\sigma Z_{0}$ where $K$ is the strike, $Z_{0}$ is the value of the underlying instrument at time $0$, $r$ is an interest rate, and $\sigma^{Ba}$ is the volatility of the Bachelier model. Then the following relation holds:
	\begin{equation}
		%\begin{cases}
		0\leq C_{S}^{Ba}(T)-C_{S}(T)\leq \E S^{\frac{3}{2}}\left(T\right)\cdot \displaystyle\frac{Z_{0}}{12\sqrt{2\pi}}\sigma^{3},\label{nierownosc1}
		%\end{cases}
	\end{equation}\label{th1}
	for the inverse subordinator $S(t)$, $t\in[0,T]$.
\end{Theorem}

\begin{Proof}
	Note that (\ref{nierownosc1}) is the generalization of the analogous relation between option prices in the classical Bachelier and B-S model \cite{Sch}:
	\begin{equation}
		0\leq C^{Ba}\left(\tau\right)-C\left(\tau\right)\leq \displaystyle\frac{Z_{0}}{12\sqrt{2\pi}}\sigma^{3}\tau^{\frac{3}{2}},\label{4}
	\end{equation}
	for $\tau\geq0.$
%=\E\left(\exp\left\{\gamma B(S(T))-\frac{\gamma^{2}}{2}S(T)\right\}(Z_{S}(T)-K)^{+}\right).\]
	Let us consider

	\begin{equation}C^{Ba}_{S}(T)=\E(C^{Ba}(S(T)))=\int_{0}^{\infty}C^{Ba}\left(x\right)g\left(x,T\right)dx,\label{niety1}
	\end{equation}
	and
	\begin{equation}C_{S}(T)=\E(C(S(T)))=\int_{0}^{\infty}C\left(x\right)g\left(x,T\right)dx,\label{niety2}\end{equation}
	where $g\left(x\right)$ - PDF of $S(T)$.
	Note that (\ref{niety1}) and (\ref{niety2}) hold by \cite{orzel} and \cite{gajda}, respectively. 
	%	where $C^{Ba}(x)$ and $C(x)$ - the fair prices of a call option in Bachelier and B-S model depending of expiration time,
	Let us observe that
	\[C_{S}^{Ba}\left(T\right)-C_{S}\left(T\right)=\int_{0}^{\infty}\left(C^{Ba}-C\right)\left(x\right)g\left(x,T\right)dx\geq0,\]
	where the equality holds by (\ref{niety1}) and (\ref{niety2}), and inequality by (\ref{4}).
	%where $\sigma^{Ba}$ is implied volatility in Bachelier Model.\newline
	Now let us use the same idea for the upper boundary of $C_{S}^{Ba}\left(T\right)-C_{S}\left(T\right)$:\newline
	\[C_{S}^{Ba}\left(T\right)-C_{S}\left(T\right)=\int_{0}^{\infty}\left(C^{Ba}-C\right)\left(x\right)g\left(x,T\right)dx \leq\int_{0}^{\infty}x^{\frac{3}{2}}g\left(x,T\right)dx\cdot \displaystyle\frac{Z_{0}}{12\sqrt{2\pi}}\sigma^{3}
	=\E S^{\frac{3}{2}}\left(T\right)\cdot \displaystyle\frac{Z_{0}}{12\sqrt{2\pi}}\sigma^{3},\]%\newline
where the first equality holds by (\ref{niety1}) and (\ref{niety2}), and the inequality by (\ref{4}). Both results can be finalized as two-sided inequalities:
	\begin{equation*}
		0\leq C_{S}^{Ba}\left(T\right)-C_{S}\left(T\right)\leq \E S^{\frac{3}{2}}\left(T\right)\cdot \displaystyle\frac{Z_{0}}{12\sqrt{2\pi}}\sigma^{3}.
	\end{equation*}

\end{Proof}\\
Note that by the relations between European put and call options for the subordinated Bachelier \cite{orzel} and the subordinated B-S model (\ref{putcallpar1}) (which for the case $r=0$ has the same form
$P_{S}=C_{S}+K-Z_{0}$), (\ref{nierownosc1}) also holds for put options. 
\begin{Example}
	Let $S(t)$ be an $\alpha$ stable inverse subordinator.
	\[ S(t)\stackrel{def}{=}S_{{\psi}_{\alpha}}\left(t\right) = \inf\{\tau \geq 0: U_{{\psi}_{\alpha}}\left(\tau\right) \geq t \}.\]
	Here $U_{{\psi}_{\alpha}}$ is the $\alpha$-stable subordinator \cite{MM, meerschaert2013inverse}.
	Note that $U_{{\psi}_{\alpha}}\left(\tau\right) < t \Longleftrightarrow S_{{\psi}_{\alpha}}\left(t\right) > \tau$. An important feature of the subordinator is that it is $1/\alpha$ self-similar, and thus its density function satisfies \cite{piryatinska2005models}
	\[ f\left(t; \tau\right) = \tau^{-1/\alpha} f\left(\displaystyle\frac{t}{\tau^{1/\alpha}}\right),\]
	where $f\left(t\right) \equiv f\left(t;1\right)$. With that, we can also obtain a similar formula for the density of the inverse stable subordinator \cite{piryatinska2005models}
	\[ g\left(\tau; t\right) = \displaystyle t^{-\alpha} g\left( \displaystyle\frac{\tau}{t^\alpha} \right), \]
	where $g\left(\tau\right) = \displaystyle\frac{1}{\alpha \tau^{1+1/\alpha}} f\left(\displaystyle\frac{1}{\tau^{1/\alpha}}\right)$. Therefore, the formula \cite{piryatinska2005models}
	\begin{equation}
		e^{-\tau s^\alpha} = \alpha \int_0^\infty e^{-st} \displaystyle\frac{\tau}{t^{\alpha+1}} g\left(\displaystyle\frac{\tau}{t^\alpha}\right) dt.
		\label{eq:laplace-subordinator}
	\end{equation}
	%%jezeli recenzent kaze rozwinac to przejscie to jest zakomentowane ponizej (to co na niebiesko)
	%\textcolor{blue}{After multiplying both sides of (\ref{eq:laplace-subordinator}) by $\tau^{k-1}$ and integrating over $\tau$, we have 
	%\begin{equation*}
	%\int_{0}^{\infty} e^{-\tau s^\alpha}\tau^{k-1} d\tau %= \alpha \int_{0}^{\infty} \int_0^\infty e^{-st} %\displaystyle\frac{\tau^{k}}{t^{\alpha+1}} %g\left(\displaystyle\frac{\tau}{t^\alpha}\right) dt %d\tau= \alpha \int_{0}^{\infty}e^{-st} \int_0^\infty % \displaystyle\frac{\tau^{k}}{t^{\alpha+1}} %g\left(\displaystyle\frac{\tau}{t^\alpha}\right) %d\tau dt=\alpha %\mathbb{E}\left(S_{{\psi}_{\alpha}}\left(1\right)^k\right) %\int_{0}^{\infty} e^{-st} \displaystyle t^{\alpha %k-1}  dt,
	%\label{eq:laplace-subordinatorextension}
	%\end{equation*}
	%where the second equality is true by Tonelli Theorem %(the funtcion $\displaystyle %e^{-st}\frac{\tau^{k}}{t^{\alpha+1}} %g\left(\displaystyle\frac{\tau}{t^\alpha}\right)$ is %non-negative and bounded for $(t,\tau)\in %[0,\infty]\times [0,\infty]$ and fixed $s>0$).
	%Thus, we have
	%\begin{equation}
	%\int_{0}^{\infty} e^{-\tau s^\alpha}\tau^{k-1} d\tau %= \alpha \mathbb{E}\left(S_{{\psi}_{\alpha}}\left(1\right)^k\rig%ht) \int_{0}^{\infty} e^{-st} \displaystyle t^{\alpha %k-1}  dt.
	%\label{eq:laplace-subordinatorextension2}
	%\end{equation}
	%After elementary transformation %(\ref{eq:laplace-subordinatorextension2}) can be %transformed into the formula for any given moment of %the inverse subordinator: %\cite{piryatinska2005models}
	%}
	Using the connection with the Laplace transform (\ref{eq:laplace-subordinator}), we can obtain the formula for any given moment \cite{piryatinska2005models}

	\begin{equation}
		\mathbb{E}\left(S(1)^k\right) = \displaystyle\frac{\Gamma\left(k+1\right)}{\Gamma\left(k\alpha+1\right)}, \quad \textrm{ for any } k > 0.
		\label{eq:moment-1-subordinator}
	\end{equation}
	Thus, by Theorem \ref{th1}, Equation (\ref{eq:moment-1-subordinator}) and by the $\alpha$ self-similarity of the inverse subordinator, we have:
	%\[ES^{\frac{3}{2}}_{\psi}\left(T\right)=\displaystyle\frac{T^{\frac{3}{2}\psi}\Gamma \left(\displaystyle\frac{3}{2}+1\right) }{\Gamma \left(\displaystyle\frac{3}{2}\psi+1\right)}\]
	\begin{equation}
		0\leq C_{S}^{Ba}\left(T\right)-C_{S}\left(T\right)\leq \displaystyle\frac{T^{\frac{3}{2}\alpha}\Gamma \left(\displaystyle\frac{3}{2}+1\right) }{\Gamma \left(\displaystyle\frac{3}{2}\alpha+1\right)}\cdot \displaystyle\frac{Z_{0}}{12\sqrt{2\pi}}\sigma^{3}.\label{nierodMM1}
	\end{equation}
	Let us observe that, for real-life parameters, the right-hand side of (\ref{nierodMM1}) is low; therefore, both models are close to each other. This result is also interesting because it indicates the relation between the solutions of two different fractional diffusion equations. Indeed, the fair price of the European call option in the subdiffusive B-S model \cite{ja2, ja3} is equal to $v(z,T)$, where $v(z,t)$ is determined by:
	\begin{equation}
		\begin{cases}
			\displaystyle {}_{0}^{c}D_{t}^{\alpha}v\left(z,t\right)=
			\frac{1}{2}\sigma^2 z^{2}\frac{\partial^2 v\left(z,t\right)}{\partial z^2}+rz\frac{\partial v\left(z,t\right)}{\partial z}-rv\left(z,t\right),\\
			v\left(z,0\right)=\max\left(z-K,0\right),\\
			v\left(0,t\right)=0,\\%\text{(bedziemy glownie rozpatrywac case D=0, D-dywidenda)},\\
			v\left(z,t\right)\sim z \quad\text{for}\quad z \rightarrow \infty, 
		\end{cases}\label{koc1}
	\end{equation}
	for $\left(z,t\right)\in(0,\infty)\times(0,T]$. Similarly, the fair price of the European call option in the subdiffusive Bachelier model is $w(z,T)$, where $w(z,t)$ is given by
	
	\begin{equation}
		\begin{cases}
			\displaystyle {}_{0}^{c}D_{t}^{\alpha}w\left(z,t\right)=
			\frac{1}{2}(\sigma^{Ba})^2 \frac{\partial^2 w\left(z,t\right)}{\partial z^2}+rz\frac{\partial w\left(z,t\right)}{\partial z}-rw\left(z,t\right),\\
			w\left(z,0\right)=\max\left(z-K,0\right),\\
			w\left(0,t\right)=0,\\%\text{(bedziemy glownie rozpatrywac case D=0, D-dywidenda)},\\
			w\left(z,t\right)\sim z \quad\text{for}\quad z \rightarrow \infty, 
		\end{cases}\label{koc2}
	\end{equation}
	for $\left(z,t\right)\in(0,\infty)\times(0,T]$. In (\ref{koc1}) and (\ref{koc2}) ''$\sim$'' denotes the limiting behavior of a function when the argument tends toward $\infty$. Furthermore, in both (\ref{koc1}) and (\ref{koc2}) ${}_{0}^{c}D_{t}^{\alpha}$ denotes a fractional Caputo derivative defined as \cite{tavares2016caputo}:
	\[{}_{0}^{c}D_{t}^{\alpha}g\left(t\right)=\frac{1}{\Gamma\left(1-\alpha\right)}\int_{0}^{t}\frac{d g\left(s\right)}{d s}\left(t-s\right)^{-\alpha}ds.\]
	Note that the interpretation of $z$ and $t$ in both (\ref{koc1}) and (\ref{koc2}) is the value of the underlying asset and the time left to $T$, respectively. The derivation of (\ref{koc2}) can be done using the same approach as in the proof of Theorem 1.3. in \cite{jatemp}. In fact, it suffices to take $\lambda=0$ and replace the standard B-S equation satisfied by $\xi(z,t)$ with the classical Bachelier PDE.
	
	Then for $r=0$, $z=K$ and $\sigma^{Ba}=\sigma z$, the difference between the solutions (\ref{koc2}) and (\ref{koc1}) for $t=T$, by (\ref{nierodMM1}) is bounded by
	\[0\leq w(z,T)- v(z,T)\leq \displaystyle\frac{T^{\frac{3}{2}\alpha}\Gamma \left(\displaystyle\frac{3}{2}+1\right) }{\Gamma \left(\displaystyle\frac{3}{2}\alpha+1\right)}\cdot \displaystyle\frac{z}{12\sqrt{2\pi}}\sigma^{3}.\]
	
\end{Example}

\section{Subordinated CRR model}
In this section, we will consider a binomial tree approximation of the subordinated B-S model introduced in Section \ref{sectionbs}.
Similarly to the classical B-S binomial tree approximation \cite{cox1979option}, we will prove that the subordinated CRR is close to the B-S model. More precisely, we will show that the underlying instrument in the subordinated CRR model converges in distribution to its equivalent in the subordinated B-S model and that in many cases the subordinated CRR option price converges to its (subordinated) B-S analogue - in both cases, the convergence is taken over the number of nodes in the binomial model. The idea presented in this paper follows the idea given by Gajda and Magdziarz article \cite{gajda}. We have that:
\begin{equation}
	\Csubbs(T) 
	= \ \E C \left(S(T)\right) = \E \big( \E^{\mathbb{Q_{B-S}}} \left( C(\tau)| \tau = S(T)\right) \big), \label{eq:exp-exp}
\end{equation}
where $\mathbb{Q_{B-S}}$ denotes the martingale measure of the classical B-S model.
We will justify that the subordinated CRR and B-S are close to each other.
%\begin{equation*}
%	C^\textrm{B-S}\left(Z_0, K, S(T), \sigma \right) \approx %\end{equation*} 
to calculate the innermost expectation on the right-hand side of formula (\ref{eq:exp-exp}), but applied to the binomial tree approximation. Then, we will use Monte Carlo simulations to calculate the sample mean using the simulated trajectories of the inverse subordinator.
An alternative view to this approach is that we are randomizing the time horizon on the classical binomial tree. Randomization is given by the inverse subordinator. Furthermore, to calculate the option price, we calculate the appropriate sample mean of prices from randomized binomial trees.

\subsection{Convergence of binomial tree asset price model to subordinated GBM}
We will consider the subordinated GBM defined in (\ref{defSubGBM}). %as
%\begin{align*}
%$\Zsub(t) = Z(S(t))$, where $Z(t)$ is the standard GBM and $S(t)$ is the inverse subordinator described in (\ref{Def_subordynatora}).
%\end{align*}
In this section, we will prove that the stock price given by the modified binomial tree converges in distribution to the subordinated GBM (conditionally on $\tau=S(t)$). %, which is given by
%\begin{align*}
%	\ln \left(\frac{\Zsub(t)}{\Zsub(0)} \right) = \ln \left(\frac{Z(S(t))}{Z(S(0))} \right) = \mu S(t) + \sigma B(S(t)),
%\end{align*}
%where $B$ is the standard Brownian motion.

Consider the realization $\tau$ of $S(t)$ for any $t\in(0,T]$. We will begin by introducing a discrete model for our approximation of the asset price. We consider the set of equidistant times $\{\tau_j = \frac{j}{n}: j =0, 1, \ldots, \floor{n \tau} \}$. Let $\ZCRR(0)$ denote the initial price of the stock $Z_{0}$, $\ZCRR(\tau_k)$ its price in ``subordinated'' time $\tau_k$. We consider a discrete conditional model of CRR type \cite{cox1979option}:
\begin{align*}
	\ln\left(\frac{\ZCRR(\tau_j)}{\ZCRR(\tau_{j-1})}\right) = \chi^{(n)}_j \in \{U,D\},
\end{align*}
where $j\in\{1, 2,\ldots, \floor{n \tau}\}$ and $\chi^{(n)}_j$ is the conditional return of the stock. Random variables $U,D$ naturally depend on $\tau$ and on the number of steps $n$ (and other parameters of the model, such as $\sigma$ and $r$).
Thus, the price of the asset $\ZCRR$ in time $\tau_k$ is given by
\begin{align*}
	\ZCRR(\tau_k) = \ZCRR(\tau_0) \prod_{j=1}^k \exp\left\{\chi^{(n)}_j\right\} \quad \textrm{ for all } k \in \{ 1, 2 \ldots, {\lfloor n \tau \rfloor}\}.
\end{align*}
To provide a probability model, we assume that $\chi^{(n)}_j$ for fixed $n$, $j= 1, 2 \ldots, {\lfloor n \tau \rfloor}$, are mutually (conditionally on $\tau=S(t)$ ) independent and identically distributed random variables in a common probability space $(\Omega, \mathcal{F}, \mathbb{P})$.%, and $n$ denotes number of steps.
That is,
\begin{align*}
	\proba\left(\chi^{(n)}_j = U\right) = p_{n} = 1 - \proba\left(\chi^{(n)}_j = D\right) \quad \textrm{ for}\quad j = 1, 2, \ldots, {\lfloor n \tau \rfloor}.
\end{align*}
\begin{Theorem}\label{frank1}
	Let us assume that the conditional expectation and variance of (conditional) returns of the stock $\chi^{(n)}_j$ are proportional to time interval $\frac{1}{n}$:
	\[  \mu_{n}=\E (\chi^{(n)}_j) = \frac{\mu}{n}, \qquad  \sigma_{n}^2=\var(\chi^{(n)}_j) = \frac{\sigma^2}{n}.\]
	Under this assumption and the previously introduced construction of a conditional binomial tree, we have \begin{equation}
	\lim_{n\to\infty}\ZCRR(S(t))\stackrel{d}{=}Z(S(t)).\label{sikort}
	\end{equation}
\end{Theorem}
\begin{Proof} For the purpose of this proof, let us define \[X_{n}(t)=\ln\left(\frac{\ZCRR(t)}{\ZCRR(0)} \right),\]
for $t\geq0$.
Let us denote: $p(x,t)$ - PDF of $X_{n}(t)$, $\vartheta\left(x,t\right)$ - PDF of $S(t)$, $h(x,t)$ - PDF of $X_{n}(S(t))$. By the total probability formula, we get:
\begin{equation} 
h(x,t)=\int_{0}^{\infty}p(x,\tau)\vartheta(\tau,t)d\tau.\label{noga1}
\end{equation} 
Let us consider the characteristic function of $X_{n}(S(t))$: \begin{multline}
\E e^{i z X_{n}(S(t))}=\int_{\mathbb{R}}e^{izx}h(x,t)dx=\int_{\mathbb{R}}e^{izx}\int_{0}^{\infty}p(x,\tau)\vartheta\left(\tau,t\right)d\tau dx=\\\int_{0}^{\infty}\int_{\mathbb{R}}e^{izx}p(x,\tau)dx\vartheta\left(\tau,t\right)d\tau=\int_{0}^{\infty}\E e^{izX_{n}(\tau)}\vartheta\left(\tau,t\right)d\tau,
\label{noga2}\end{multline}	
where the second equality holds by  (\ref{noga1}). Furthermore, we could change the integration order in the third equality of (\ref{noga2}) by the Fubini Theorem, because the densities are non-negative a.s. and integrable, and the absolute value of the characteristic function is bounded by $1$. By convergence of the underlying instrument in the standard CRR model \cite{cox1979option} we get \[\lim_{n\to\infty}X_{n}(t)\stackrel{d}{=}\tilde{B}(t),\]
where $\tilde{B}(t)=\mu + \sigma \sqrt{t} Z$ and $Z \sim \mathcal{N}(0,1)$. To finish the proof, it suffices to show the convergence of the characteristic functions of $X_{n}(S(t))$ to the characteristic function of $\tilde{B}(S(t))$:
\begin{equation*}
\lim_{n\to\infty}\E e^{izX_{n}(S(t))}=\lim_{n\to\infty}\int_{0}^{\infty}\E e^{izX_{n}(\tau)}\vartheta\left(\tau,t\right)d\tau=\\\int_{0}^{\infty}\lim_{n\to\infty}\E e^{izX_{n}(\tau)}\vartheta\left(\tau,t\right)d\tau=\int_{0}^{\infty}\E e^{iz\tilde{B}(t)}\vartheta\left(\tau,t\right)d\tau=\E e^{iz\tilde{B}(S(t))},
\end{equation*}	
where the second equality holds by Lebesgue's Dominated Convergence Theorem.
\end{Proof}

Theorem \ref{frank1} can be strengthened for the case of functional convergence in the Skorokhod space. To formulate the corresponding theorem, let us first introduce the Definitions \ref{skd1}-\ref{skd3} \cite{whitt2002stochastic}.
\begin{Definition}\label{skd1}
	The Skorokhod space $D^{k}\vcentcolon= D(\mathbb{I},\mathbb{R}^{k})$ is a space of all functions of $c\grave{a}dl\grave{a}g$-type defined on some interval $\mathbb{I}\subset\mathbb{R}$ having values in $\mathbb{R}^{k}$. For $k=1$, we use the notation $D=D^{1}$.
\end{Definition}
\begin{Definition}\label{skd2}
	The metrics $M_{1}$ and $M_{2}$ we define by
	\[d_{M_{1}}(f_{1},f_{2})\vcentcolon=\inf_{\substack{(u_{j},r_{j})\in \pi(f_{j})\\ j=1,2}}\left(\max\left(\norm{u_{1}-u_{2}}_{\infty},\norm{r_{1}-r_{2}}_{\infty}\right)\right)\]
	\[d_{M_{2}}(f_{1},f_{2})\vcentcolon=m_{H}(\Gamma_{f_{j}},\Gamma_{f_{2}}),\]
	where $\pi(f)$ is a set of parametric representations of $f$ in $D$ and $\Gamma_{f}$ is a completed graph defined as:
	\[\Gamma_{f}\vcentcolon=\{(z,t)\in R\times[0,1]:z=\alpha f(t-)+(1-\alpha)f(t)\},\] for some $\alpha\in[0,1]$ and $f\in D([0,1],R)$. Moreover, $m_{H}$ is a Hausdorff's metrics defined as \[m_{H}(K_{1},K_{2})\vcentcolon=\max\left(\sup_{x_{1}\in K_{1}}m(x_{1},K_{1}),\sup_{x_{2}\in K_{2}}(x_{2},K_{2})\right),\] for the distance of point $x$ from set $A$ defined by \[m(x,A)\vcentcolon=\inf_{y\in A}m(x,y),\]
	where $m(x,y)$ is a standard metrics in $\mathbb{R}^{k}$.\label{buza}
\end{Definition}
\begin{Definition}\label{skd3}
	Let $D_{0}$ be a subset of all $f\in D^{k}$, where for each $i=1\ldots,k$ $f^{(i)}\geq 0$, where $f^{(i)}$ is the $i$-th coordinate of $f$. Then $D_{\uparrow}$ denotes the set of functions belonging to $D_{0}$, which do not decrease with respect to each coordinate.
\end{Definition}
	Let us present a Theorem $13.2.4$ from \cite{whitt2002stochastic}, which will be used to prove a Theorem \ref{frank11}:
	\begin{Theorem}
		Suppose that $(x_{n},y_{n})\to (x,y)$ in $D^{k}\times D^{1}_{\uparrow}$. If\\
	%	\begin{enumerate}
			%\item 
			\indent 1. $y$ is continuous and strictly increasing at $t$ whenever $y(t)\in Disc(x)$ and\\
			%\item
			\indent 2. $x$ is monotone on $[y(t-),y(t)]$ and $y(t-)$, $y(t)]\notin Disc(x)$ whenever $t\in Disc(y)$,\\
		%\end{enumerate} 
		then $x_{n}\circ y_{n}\to x\circ y$ in $D^{k}$, where the topology throughout is $M_{1}$ or $M_{2}$.\label{TWduka}		
	\end{Theorem}
In Theorem \ref{TWduka}, $Disc(f)$ denotes a set of discontinuities of $f$.
\begin{Theorem}\label{frank11}
Under the assumptions of Theorem \ref{frank1} there holds a functional convergence in Skorokhod space $D([0,T],\mathbb{R})$ of the price of the underlying instrument in the subordinated CRR model to the price of this asset in subordinated B-S model: \begin{equation*}
		\lim_{n\to\infty}\ZCRR(S(t))\stackrel{D}{=}Z(S(t)),\label{sikort1}
		\end{equation*}
		 where the topology throughout is $M_{1}$ or $M_{2}$ introduced in Definition \ref{buza}.
\end{Theorem}
\begin{Proof}
	Let us denote $x_{n}=\ln\left( \frac{\ZCRR(\cdot)}{\ZCRR(0)}\right)$, $y_{n}=S(\cdot)$, $x(t)=t \mu + \sigma B(t)$, $y(t)=S(t)$, where $B(t)$ is a Brownian motion. Observe that due to the convergence of the price of the underlying instrument in the classical CRR model to the price of this asset in the standard B-S model in $D([0,T],\mathbb{R})$ \cite{di2011advanced}, we have $\lim_{n\to\infty}x_{n}\stackrel{D}{=} x$. All other assumptions of Theorem \ref{TWduka} are satisfied because $B$ and $S$ have continuous trajectories and $S$ is nondecreasing. So, based on Theorem \ref{TWduka}, the proof is completed. Let us observe that the limit process has a continuous distribution, hence there holds also a convergence of finite-dimensional distributions.
	\end{Proof}

\subsection{Convergence of the  subordinated binomial tree model to subordinated B-S model}
 % - the single-period interest rate.%, and $r$ - the continuously compounded risk-free interest rate.
We will use the following notations: $q$ - risk neutral probability of up-movement on a binomial tree, $\mathbb{H}$ - a risk neutral measure corresponding to $q$ and $1-q$, $U$ and $D$ - parameters of the binomial tree model (dependent on $n$ and $\tau=S(T)$).	By the basic properties of the classical binomial model, we get \cite{cox1979option}: 
 		 \begin{equation}
 			q=\frac{R-D}{U-D},\label{www0}
 		\end{equation}
 	\begin{equation}
 		U=e^{\sigma \sqrt{\frac{\tau}{n}}},\label{www1}
 	\end{equation}	
 	\begin{equation}
 		D=e^{-\sigma \sqrt{\frac{\tau}{n}}},\label{www2}
 	\end{equation}
 	where $\tau$ is a realization of $S(T)$.	Since we consider a discrete model, we assume that
 	\begin{equation}
 		R=(1+r)^{1/n}.\label{www3}
 	\end{equation}
 
 Let us focus on the European option with the payoff function $f(Z(T))$ (where $Z(T)$ is the price of the stock at time $T$) in the subordinated CRR and subordinated B-S model. Then, we have
\begin{Theorem}
	\label{proba:2}
	%and $\sigma(S(\tau))_{\tau\in(0,T]}$ - $\sigma$-algebra generated by $S(\tau)$ for $\tau\in(0,T]$. 

	Let us assume that there exists a constant $c>0$ such that for $x\in\mathbb{R}$ the payoff function $f(x)$ follows: 
	\begin{equation}
	    f(x)\leq c\left|x \right|.\label{tire1}
	\end{equation}
	Then, according to the assumptions of Theorem \ref{frank1}, the price of the European option in the subordinated binomial model, which is given by the following formula:
	\begin{equation*}
		\Csubbin (T)= \E\big(\E^{\mathbb{H}}\big( \frac{1}{R^{n}}\sum_{j=0}^{n} {{n} \choose j} q^j (1-q)^{{n}-j}f(\ZCRR(0)U^{j}D^{n-j})\big|\tau=\\S(T) \land \ZCRR(0)=Z_{0}\big)\big),
	\end{equation*}
	converges to the option price in the subordinated B-S model: %\cite{gajda}
	\begin{equation*}
		\Csubbs (T) = \E \left(\E^{\mathbb{Q_{B-S}}}\left(e^{-r\tau} f(Z(\tau))|\tau=S(T) \land Z(0)=Z_{0}\right) \right),
	\end{equation*}
	where $\mathbb{Q_{B-S}}$ denotes the martingale measure of the classical B-S model.
\end{Theorem}
\begin{prof}
	By the convergence of the CRR model to the B-S  \cite{cox1979option}, for expiration time $\tau$ we have:
	\begin{equation}
		\lim_{n\to\infty}\E^{\mathbb{H}} \left(\frac{1}{R^{n}}\sum_{j=0}^{n} {{n} \choose j} q^j (1-q)^{{n}-j}f(\ZCRR(0)U^{j}D^{n-j}) \big| \ZCRR(0)=Z_{0}\right)= e^{-r\tau}\E^{\mathbb{Q_{B-S}}}\left( f(Z(\tau)| Z(0)=Z_{0} \right).\label{conditiontau}
	\end{equation}
	Therefore, using the fact that for $t\in(0,T]$ there holds (\ref{sikort}), %\textcolor{blue}{  \[\lim_{n\to\infty} \ZCRR(S(t))\stackrel{d}{=} Z(S(t)),\]} 
	after conditioning (\ref{conditiontau}) by $\tau=S(T)$ we have 
	\begin{equation}
		\lim_{n\to\infty}\E^{\mathbb{H}} \big(\frac{1}{R^{n}}\sum_{j=0}^{n} {{n} \choose j} q^j (1-q)^{{n}-j}f(\ZCRR(0)U^{j}D^{n-j})\big|\tau=S(T) \land \ZCRR(0)=Z_{0}\big)\stackrel{d}{=} \E^{\mathbb{Q_{B-S}}}\big(e^{-r\tau} f(Z(\tau)|\tau=S(T)\land Z(0)=Z_{0} \big).\label{flir3}
	\end{equation}
	After taking the expected value on both sides of (\ref{flir3}), we get
	\begin{equation}
		\E(\lim_{n\to\infty}\E^{\mathbb{H}} \big(\frac{1}{R^{n}}\sum_{j=0}^{n} {{n} \choose j} q^j (1-q)^{{n}-j}f(\ZCRR(0)U^{j}D^{n-j})\big|\tau=S(T) \land \ZCRR(0)=Z_{0}\big))=\\ \E(\E^{\mathbb{Q_{B-S}}}\big(e^{-r\tau} f(Z(\tau)|\tau=S(T)\land Z(0)=Z_{0} \big)).\label{flir5}
	\end{equation}
	To finish the proof, we have to show that  \begin{multline}
	\E(\lim_{n\to\infty}\E^{\mathbb{H}} \big(\frac{1}{R^{n}}\sum_{j=0}^{n} {{n} \choose j} q^j (1-q)^{{n}-j}f(\ZCRR(0)U^{j}D^{n-j})\big|\tau=S(T) \land \ZCRR(0)=Z_{0}\big))=\\
	\lim_{n\to\infty}\E(\E^{\mathbb{H}} \big(\frac{1}{R^{n}}\sum_{j=0}^{n} {{n} \choose j} q^j (1-q)^{{n}-j}f(\ZCRR(0)U^{j}D^{n-j})\big|\tau=S(T) \land \ZCRR(0)=Z_{0}\big)).\label{hecer}\end{multline}
	By (\ref{tire1}), we have
	\begin{multline}
	    \left|\frac{1}{R^{n}}\sum_{j=0}^{n} {{n} \choose j} q^j (1-q)^{{n}-j}f(\ZCRR(0)U^{j}D^{n-j}\right|\leq   \frac{1}{R^{n}}\sum_{j=0}^{n} {{n} \choose j} q^j (1-q)^{{n}-j}c\left|\ZCRR(0)U^{j}D^{n-j}\right|=\\\frac{c\left|\ZCRR(0)\right|}{R^{n}}\sum_{j=0}^{n} {{n} \choose j} (q\left|U \right|)^{j}((1-q)\left|D\right|)^{n-j}=\frac{c\left|\ZCRR(0)\right|}{R^{n}}((q\left|U \right|)+((1-q)\left|D\right|))^{n}\\=\frac{c\left|\ZCRR(0)\right|}{R^{n}}\left(\frac{R-U^{-1}}{U-U^{-1}}U +\frac{U-R}{U-U^{-1}}U^{-1}\right)^{n}=c|\ZCRR(0)|,\label{fiubz}
		\end{multline}
	where the third equality is satisfied by (\ref{www0})-(\ref{www2}),and the last equality is satisfied by (\ref{www3}).
Since we bounded\\ $\displaystyle\left|\frac{1}{R^{n}}\sum_{j=0}^{n} {{n} \choose j} q^j (1-q)^{{n}-j}f(\ZCRR(0)U^{j}D^{n-j}\right|$ by the integrable function independent of $n$, by the dominated convergence theorem, (\ref{hecer}) holds. Therefore, the proof is completed.
\end{prof}

Analogously, the convergence of prices of a wide range of options can be shown. Indeed, it is enough to change in the proof the option pricing formula of the CRR model. Note that (\ref{tire1}) is satisfied, for example, for payoffs of the vanilla call (or put) options. However, (\ref{tire1}) does not hold, e.g., for the lookback option considered in Figure \ref{figcrr4}. Figure\mbox{ }\ref{figcrr4} indicates that the output of the CRR method is close to the real value. However, the use of the binomial method for such options should be avoided as long as the theoretical justification is missing.

Thus, we have the following algorithm for pricing any option in the subordinated B-S model:%\newpage 

\begin{algorithm}[h!]
	\caption{}%\label{euclid}
	\begin{algorithmic}[1]
		\For{$i$ from $1$ to $M$}
		\State Compute $S^{(i)}(T)$, $PS^{(i)}=P(S^{(i)}(T))$,\\ where $P(T)$ - value of an option with expiration time $T$ in the classical CRR model, $S^{(i)}(T)$ - the $i$-th realization of inverse subordinator at $T$
		\EndFor\label{euclidendwhile}
		\State $Price = mean([PS^{(1)},\ldots,PS^{(M)}])$ 
		\State \textbf{Return} $Price$.
	\end{algorithmic}
\end{algorithm}

We recall the observation of \cite{heston2000rate} that the convergence rate depends on the smoothness of the option payoff
functions. Since most of the payoff functions are not continuously differentiable (e.g., in the case of European options the critical point is at the strike), the rate of convergence can be much lower than is assumed by many academics and practitioners. This problem can resolve the smoothing of the payoff function, for example, by transforming the original payoff function $f(x)$ into $\displaystyle f^{\ast}(x)=\frac{1}{2\Delta x}\int_{-\Delta x}^{\Delta x} f(x-y)dy$ \cite{heston2000rate} for space step $\Delta x$. The same remark is true for both the CRR method and the finite difference (FD) method \cite{ja3, ja2} (the reference method for numerical examples). The impact of smoothing the payoff function on numerical methods will not be further investigated in this paper, and therefore the smoothing techniques will not be used in numerical examples.

\section{Numerical examples}

\begin{figure}[ht]
	%\centering
	\raggedleft
	\includegraphics[scale=0.36]{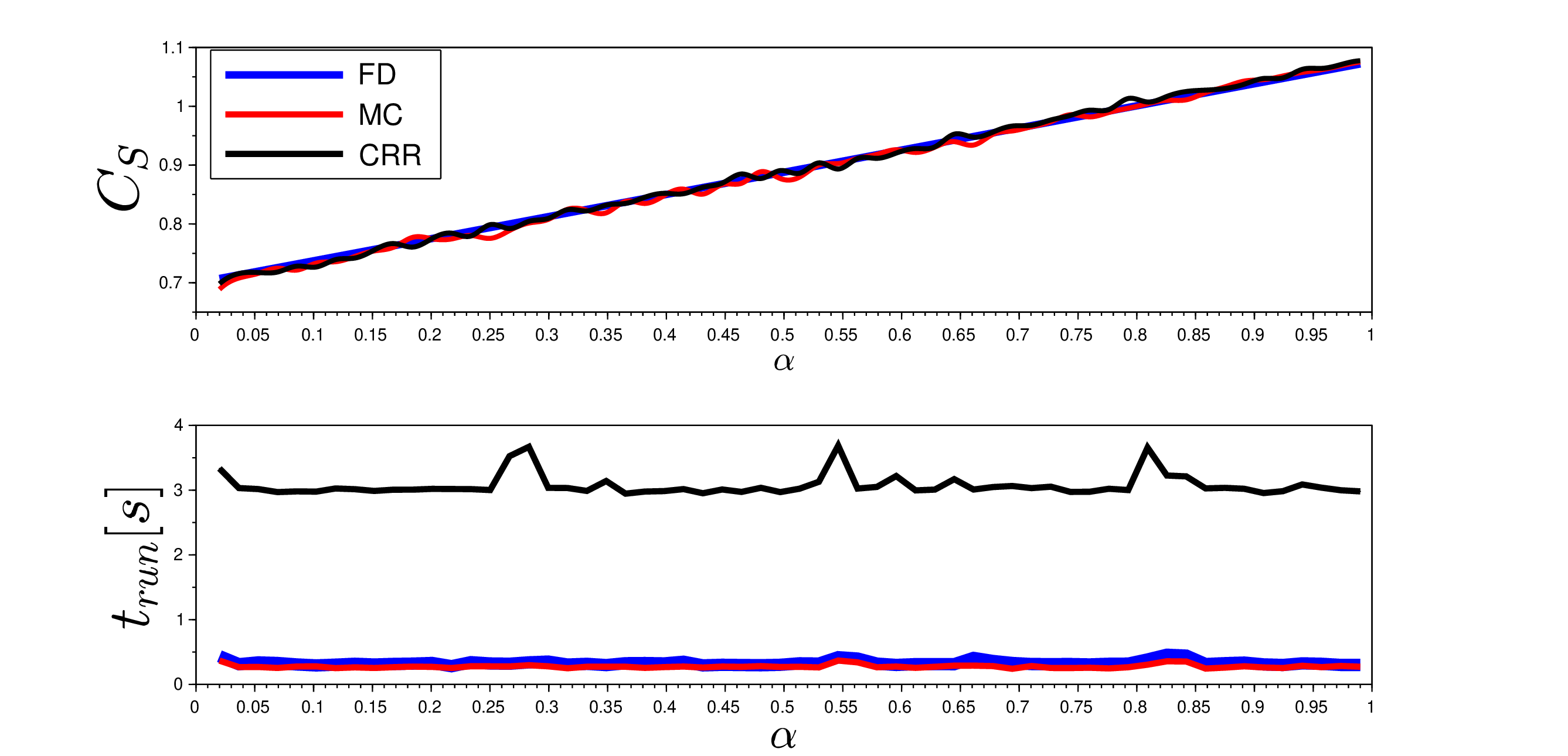}
	\caption{{\bf The European call option in subdiffusive B-S computed by different methods (top panel) with the corresponding time of computation (bottom panel).} 
	}  \label{figcrr2} 
\end{figure}

\begin{figure}[ht!]
	\centering
	\raggedleft
	\includegraphics[scale=0.36]{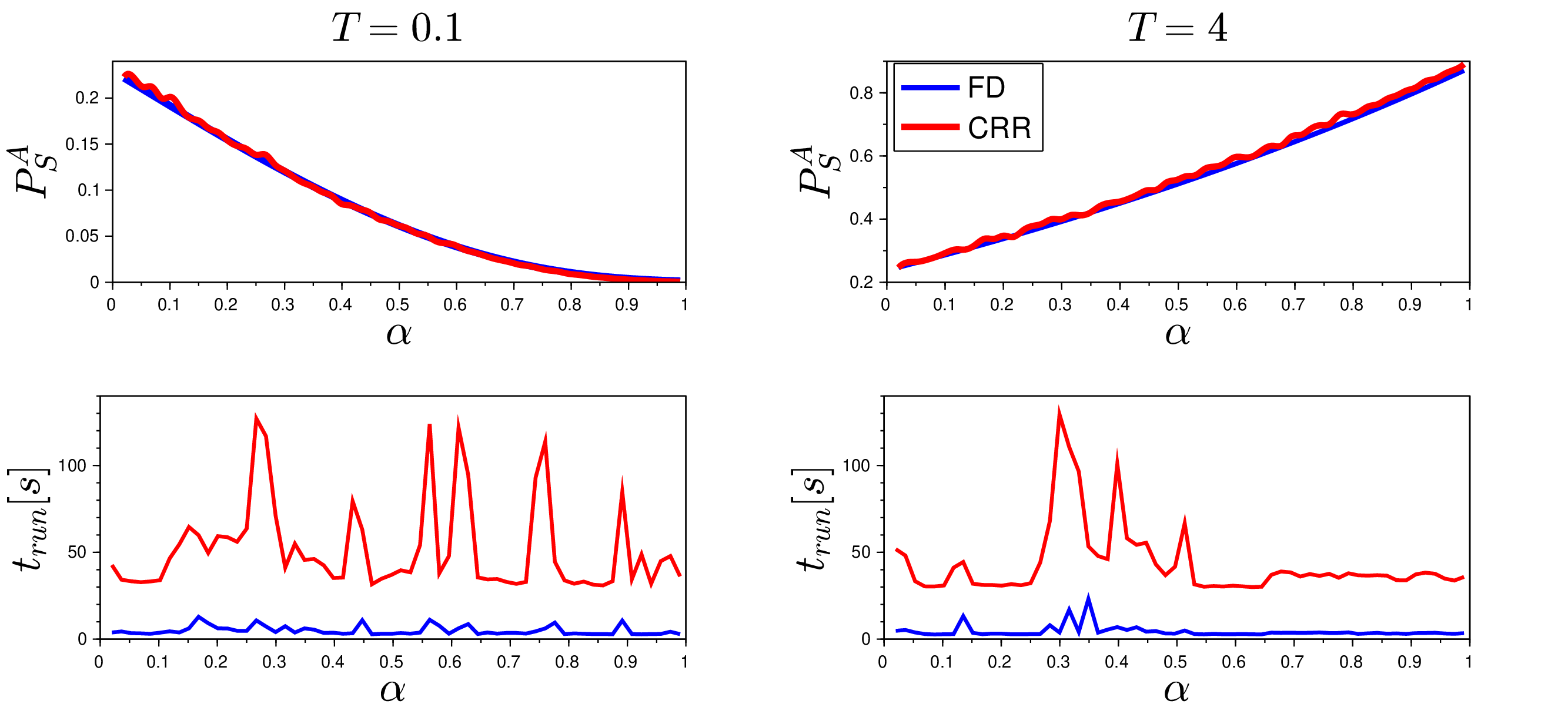}
	\caption{{\bf The American put option in subdiffusive B-S model computed by different methods (top panels) with the corresponding time of computation (bottom panels).}}   \label{figcrr3} 
\end{figure}

\begin{figure}[ht]
	%\centering
	\raggedleft
	\includegraphics[scale=0.39]{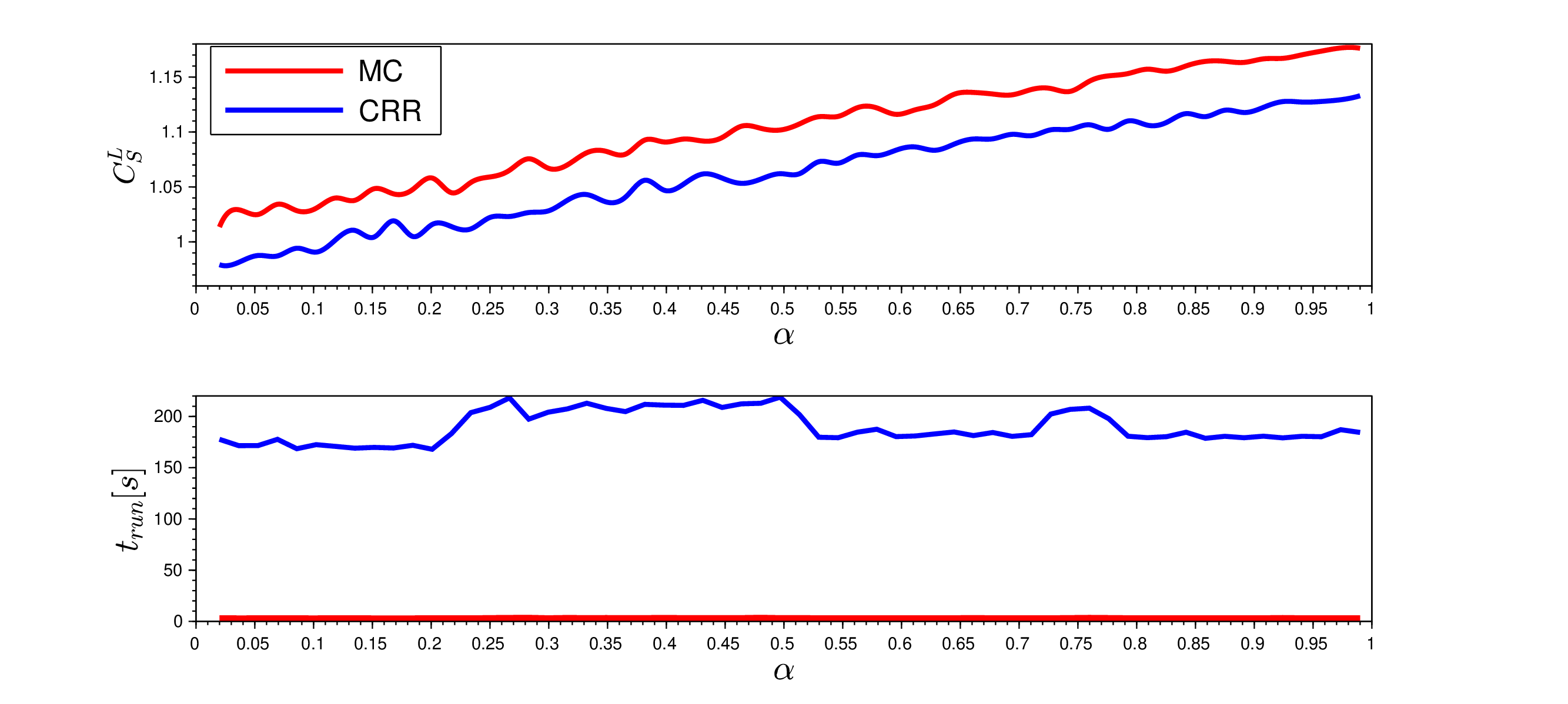}
	\caption{{\bf The European lookback call option in subdiffusive B-S computed by different methods (top panel) with the corresponding time of computation (bottom panel). }} \label{figcrr4}   
\end{figure}

Let us focus on the $\alpha$-stable inverse subordinator. We compare the numerical methods for computing the prices of European call ($C_{S}$), American put ($P^{A}_{S}$), and European lookback call option ($C^{L}_{S}$) with floating strike (with the payoff function $Z_{\psi_{\alpha}}(T)-\min_{t\in[0,T]}(Z_{\psi_{\alpha}}(t)$ ) with regard to parameter $\alpha$ for Figures \ref{figcrr2} -\ref{figcrr4} respectively. We take $r=0.04$, $\sigma=1$ and for the FD method $x_{max}=10$, $x_{min}=-20$, $\theta=0$ \cite{ja3, ja2}.
For vanilla options, we choose as the reference FD method \cite{ja3,ja2}, for the lookback option the Monte Carlo (MC) method explained in Appendix B. In Figure \ref{figcrr2} a MC method introduced in \cite{MM} is also presented. 

In Figure \ref{figcrr2} all methods converge to the same result, but CRR is slower than the FD and MC methods.
The parameters are $T=2$, $Z_{0}=2$, $K=2$. For FD, the number of time and space nodes are $N=120$, $m=80$, respectively. For CRR, the number of tree nodes $n=100$. 
For MC and CRR, $M=3000$. In Figure \ref{figcrr3} we compare CRR with the FD method introduced in \cite{ja3}. Both methods match the real solution, but the CRR is significantly slower than the FD.  It is worth mentioning that the CRR works clearly better than the Longstaff-Schwarz (L-S) method introduced in \cite{longstaff2001valuing}.  In contrast to L-S, the CRR method also works correctly for small values of $T$ and $\alpha$ (see Figure $7$ in \cite{ja3}). 
The parameters are $Z_{0}=5$, $K=2$, $N=170$, $m=200$, $n=100$, $M=3000$. To prepare the Figure \ref{figcrr4}, for each realization of $S_{{\psi}_{\alpha}}(t)$ we are looking for the minimum iteratively, as is proposed in \cite{cheuk}. The parameters are $Z_{0}=2$, $n=80$,
$T=1$ and for both methods $M=7000$. In Figures \ref{figcrr2}, \ref{figcrr3} and \ref{figcrr4} we see that CRR is the slowest/least precise presented numerical method (increasing $n$ and $M$ in each case follows getting closer to the real solution, but also increasing the time of computation). Note that the greatest advantage of the introduced binomial approach is that it can be applied to a wide range of options and inverse subordinators. 
\section*{Conclusions}
In this paper, we have introduced the subordinated CRR model. We have shown that in terms of convergence in the distribution and in Skorokhod space, of the underlying assets and convergence of many option prices (in both cases with respect to the number of nodes), this model is close to the subordinated B-S model. This is a generalization of the celebrated Cox-Ross-Rubinstein result \cite{cox1979option} and gives motivation to approximate an option price in the subordinated B-S by the binomial MC method proposed in this paper. Furthermore, we have generalized the Schachermayer-Teichmann result \cite{Sch} by providing the analytical relation between the prices of European options in the subordinated Bachelier and B-S models. As a supplement to the main course of the paper, we have also derived the PDE describing the price of the European call option in the subdiffusive Bachelier model, and we have presented the Monte Carlo method for valuing lookback options in the subordinated B-S model. Finally, we present some numerical examples.
%\appendix

%\section{Monte Carlo for lookback option}
\paragraph*{Appendix A}

\begin{Theorem}
	Let $C^{L}(T)$ be a fair price of the European lookback call option with floating strike in the B-S model, given by the semidirect formula \cite{heuwelyckx2014convergence}:
	\begin{equation*}
		\displaystyle C^{L}(T)= Z_{0}\left(1+\frac{\sigma^2}{2r}\right)\Phi(a_{1})-Z_{0}e^{-rT}\left(1-\frac{\sigma^2}{2r}\right)\Phi(a_{2})-Z_{0}\frac{\sigma^2}{2r},
	\end{equation*}
	where $a_{1}=\left(r/\sigma+\sigma/2\right)\sqrt{T}$, $a_{2}=\left(r/\sigma-\sigma/2\right)\sqrt{T}$ and $\Phi(x)$ is the CDF of the normal distribution. Then its subordinated equivalent $C_{S}^{L}(T)$ is equal: 
	\begin{equation*}
		C_{S}^{L}(T)=\E C^{L}(S(T)).
	\end{equation*}
	Moreover, if we consider the $\alpha$-stable inverse subordinator $S(t)=S_{{\psi}_{\alpha}}(t)$, then
	\begin{equation*}
		C_{S}^{L}(Z_{0},r,\sigma,T)=\int_{0}^{\infty}C^{L}(x)T^{-\alpha}g_{\alpha}(x/T^{\alpha})dx,
	\end{equation*}
	where $g_{\alpha}(z)$ is given in terms of the Fox function $H^{1,0}_{1,1}\big (z^{(1-\alpha,\alpha)}_{(0,1)}\big)$ (see \cite{MM} and the references therein).
\end{Theorem}

\begin{Proof}
	The proof is analogous to the proof of Theorem $4$ from \cite{MM}.
\end{Proof}\\
Note that the analogical result can be obtained for other lookback options. The first statement provides an easy and fast numerical method, we only have to approximate the expected value by its finite equivalent. Therefore, we have 
\begin{equation*}
	C_{S}^{L}(T)\approx\frac{\displaystyle\sum_{i=1}^{M}C^{L}(S^{(i)}(T))}{M},
\end{equation*}
where $S^{(i)}(T)$ ($i=1,\ldots,M$) are independent realizations of $S(T)$ and ``$\approx$'' denotes an approximation (which reduces to the equality for $M\to\infty$). Note that the Monte Carlo methods considered in this paper can be improved, e.g., using variance reduction techniques \cite{zhang2020value}.

% Either type in your references using
% \begin{thebibliography}{}
% \bibitem{}
% Text
% \end{thebibliography}
%
% or
%
% Compile your BiBTeX database using our plos2015.bst
% style file and paste the contents of your .bbl file
% here. See http://journals.plos.org/plosone/s/latex for 
% step-by-step instructions.
%
%\section*{Author Contributions}
%Micha\l{} Balcerek:  Validation, Formal analysis, Resources, Methodology,  Investigation, Writing- Original draft preparation.

%Grzegorz Krzy\.zanowski: Resources, Conceptualization, Methodology, Software, Validation, Formal analysis, Investigation, Visualization, Writing- Original draft preparation.

%Marcin Magdziarz: Conceptualization, Supervision, Formal analysis, Methodology, Writing - Review and Editing.
\section*{Acknowledgements}
The research of G.K. and M.M. was partially supported by the NCN Sonata Bis 9 grant nr 2019/34/E/ST1/00360.

%\bibliographystyle{plain}
%\bibliography{main.bib}

\end{document}